\documentclass[%
reprint,
amsmath,amssymb,
aps,
longbibliography,
onecolumn
]{revtex4-1}

\usepackage[resetlabels]{multibib}
\usepackage{multirow}
\newcites{sec}{Reference}

\usepackage{graphicx}
\usepackage{dcolumn}
\usepackage{hyperref}
\hypersetup{colorlinks=true, citecolor=blue, urlcolor=blue, linkcolor=blue}
\usepackage{bm}
\usepackage{epstopdf}
\usepackage{titlesec}
\titleformat{\chapter}[display]
{\normalfont\Large\bfseries}{\thechapter}{11pt}{\Large}
\titleformat{\section}
{\normalfont\large\bfseries}{\thesection}{11pt}{\large}
\titlespacing*{\chapter}{0pt}{0pt}{15pt} 
\titlespacing*{\section}{0pt}{3.5ex plus 1ex minus .2ex}{2.3ex plus .2ex}
\usepackage[titletoc]{appendix}
\usepackage{color}
\usepackage{algorithm}
\usepackage{algpseudocode}

\usepackage{enumitem}

\floatname{algorithm}{Algorithm}

\begin{document}

\preprint{APS/123-QED}

\title{Signal Lasso with Non-Convex Penalties for Efficient Network Reconstruction and Topology Inference}

\author{Lei Shi$^{1,2}$}
\email{lshi@ynufe.edu.cn}
\author{Jie Hu$^3$}
\author{Huaiyu Tan$^1$}
\author{Libin Jin$^2$}
\author{Wei Zhong$^4$}
\author{Chen Shen$^{5}$}

\affiliation{
	1. School of Statistics and Mathematics, Yunnan University of Finance and Economics, Kunming, 650221, China. \\
	2. Interdisciplinary Research Institute of Data Science, Shanghai Lixin University of Accounting and Finance, Shanghai 201209, China.\\
	3. Management School, Science and Technology University of China,Hefei, China\\
	4. Department of Statistics and Data Science, School of Economics, Xiamen University,  Xiamen, 361005, China.\\
	5. Faculty of Engineering Sciences, Kyushu University, Kasuga-koen, Kasuga-shi, Fukuoka 816-8580, Japan}
\date{\today}

\begin{abstract}
Inferring network structures remains an interesting question for its importance on the understanding and controlling collective dynamics of complex systems. The most of real networks exhibit sparsely connected properties, and the connection parameter is a signal connected (0 or 1). The existing shrinking methods such as Lasso-type estimation can not suitably reveal such property.  A new method, called by {\it signal lasso} (or its updating version: adaptive signal lasso) was proposed recently to solve the network reconstruction problem, where the signal parameter can be shrunk to either 0 or 1 in two different directions,  and they found signal lasso outperformed lasso-type method. The signal lasso or adaptive signal lasso employed the additive penalty of signal and non-signal terms which is a convex function and easily to complementation in computation. However their methods need tuning the one or two parameters to find an optimal solution, which is time cost for large size network. In this paper we propose new signal lasso method based on two penalty functions to estimate the signal parameter and uncovering network topology in complex network with a small amount of observations. The penalty functions we introduced are non-convex function, thus  coordinate descent algorithms are suggested. We find in this method the tuning parameter can be set to a large enough values such that the signal parameter can be completely shrunk either 0 or 1.  The extensive simulations are conducted in linear regression models with different assumptions, the evolutionary-game-based dynamic model and Kuramoto model of synchronization problem. The advantage and disadvantage of each method are fully discussed in various conditions. Finally a real example comes from behavioral experiment is used for illustration.  
Our results show that signal lasso with non-convex penalties is effective and fast  in estimating signal parameters in linear regression model.

\end{abstract}

\maketitle
\section{Introduction}

Complex network has wide application in many fields and has been made great progress recently (\cite{r2}; \cite{r3}). In a complex network, the pattern of the node-to-node interaction or the network topology is unknown, and uncovering the network topology based on a series of observable quantities obtained from experiments or observations is very important and plays a potential role for understanding and controlling collective dynamics of complex systems (\cite{r14}; \cite{r26}; \cite{r24}). Network reconstruction is an inverse problem in network science and has received great attention recently. Some typical examples include the reconstruction of gene networks using expression data in biology (\cite{r12}; \cite{r13}), extraction of various functional networks in the human brain from activation data in neuroscience (\cite{r31}; \cite{r32}; \cite{r44}), and detection of  organizational networks in social science and trade networks in economics. Recently evolutionary-game based dynamics has also been used to study the network reconstruction, in which it is possible to observe a series of discrete but a small amount of quantities (\cite{r26}; \cite{r14};  \cite{r38, r39}). In this case, the problem can be transformed into a statistical linear model with sparse and high-dimensional properties.

To understand how such signal parameter appeared in real practice, we use two typical  examples to illustrate the problem. The first is a dynamic equation governed the evolution state in general complex systems, which can be written as following differential equations (\cite{r27}; \cite{r35})
\begin{equation}\label{eq0}
\dot{\bf y}_i(t)=\psi_{i}({\bf y}_i(t),\nu_i)+\sum_{i=1}^N a_{ij}\phi_{ij}({\bf y}_i(t),{\bf y}_j(t))+\epsilon_i(t),
\end{equation}
where ${\bf y}_i(t)$ denote a $m$-dimensional internal stats variable of a system consisting of $N$ dynamic units at time $t$, the function $\psi_{i}\in R^m$
and $\phi_{ij}$  $\in R^m$ respectively define the intrinsic and interaction dynamics of the units. $\epsilon_i(t)$ is a dynamic noise term and 
$\nu_i$ is a set of dynamic parameters. Finally, the term $a_{ij}$ defines the interaction topology and is called by adjacency matrix such that $a_{ij}$ = 1 if there is a direct physical interaction from unit $j$ to $i$ and $a_{ij}$ = 0 otherwise. The matrix $A=[a_{ij}]$ completely defines a network with size $N$, that is, an abstraction used to model a system that contains discrete interconnected elements. The elements are represented by nodes (also called vertices) and connections are represented by edges. In general ${\bf y}_i(t)$ can be observed as time series data but $a_{ij}$ for $i=1, \cdots, N$ are unknown and need to be estimated. It is clear Eq. (\ref{eq0}) can be rewritten as linear regression model if functional form of $\phi_{ij}$ and $y_{i}(t)$ are known. This model include well known synchronization model, oscillator networks, spreading network, and so on (see \cite{r27}, \cite{r39, r40}).

Another example is the network reconstruction based on evolutionary game on structured populations, where node represent a player, link means the corresponding two players have game relationship. In this model, the payoff and strategy of each player taken from interaction with other players can be observed but the relationship between players are unknown and need to be estimated (\cite{r26}; \cite{r14};  \cite{r24}; \cite{r39}). 
We use a prisoner's dilemma game (PDG) as an example with the payoff matrix defined by
\begin{equation}\label{eq1}
\mathbf{M_{PDG}} =
\left( \begin{array}{cc}
R \ \ S \\
T \ \ P
\end{array}\right),
\end{equation}
This game is characterized with the temptation to defect $T$, reward for mutual cooperation $R$, punishment for mutual defection $P$, as well as the sucker's payoff $S$.  These quantities satisfies that $T>R>P>S$, and the mutual defection are the equilibrium solutions of PDG (\cite{r36}; \cite{r23}; \cite{r17}; \cite{r42}).  In some experiment driving studies, each player can interact with other players by choosing either a cooperator (C) or defector (D) to obtain their payoff and procedure is continued until to some predetermined number (\cite{r19}; \cite{r33, r34}). Let $s_i=(1,0)$ indicates that player $i$ takes action C ( or $s_i=(0,1)$ means taking action D). In spatial PDG game, the $i$ player, say focal player, acquires its fitness (total payoff) $F_i$ by playing the game with all its connected neighbors, which is defined by:
\begin{equation}\label{eq3}
F_i=\sum_{j \in \Omega_i} s_{i}M_{PDG}s_{j}=\sum_{j=1, j \neq i}^N a_{ij} P_{ij},
\end{equation}
where $\Omega_i$ represents the set of all connected neighbors of player $i$, $P_{ij}=s_{i}M_{PDG}s_{j}$.  
Eq. (\ref{eq3}) can be converted into a linear model, where $a_{ij}$ is the elements of adjacency matrix for a given network. $a_{ij}=1$ means player $i$ and $j$ are connected, otherwise ($a_{ij}=0$) they are not connected. When the process continue forward to produce
a time series data, we obtained a linear model with signal parameters (0 or 1) as regression coefficients. The detailed description of this framework will be presented by a real example of human behavioral experiments in Section \ref{sec5} of this paper.

In the fields of complex systems and applied physics as illustrated by above two examples, traditional estimation methods such as compressed sensing (CS) and the lasso method have been widely used for network reconstruction (\cite{r26}; \cite{r14}; \cite{r38}). These studies found that the lasso method is particularly robust against noise in the data, making it effective for reconstructing sparse signals. Given that the true value of $\hat{a}_{ij}$ should be either 0 or 1, they determined that player $i$ has a game relationship (connection) with player $j$ if $\vert \hat{a}_{ij} - 1 \vert \leq 0.1$, and no relationship if $\vert \hat{a}_{ij} \vert \leq 0.1$. Here, $\hat{a}_{ij}$ is the estimator of $a_{ij}$. In cases where the estimator falls outside these ranges, the relationship between the players is not identifiable.
Although CS and the lasso method are effective at shrinking parameter estimates toward zero, especially in sparse complex networks, they often fail to shrink existing links to their true value of 1, which can reduce estimation accuracy. To address this limitation, \cite{r39} introduced the signal lasso (SigL) method, which adds a penalty for the signal parameter. They found that SigL outperformed both lasso and CS in network reconstruction. \cite{r21} applied the SigL method to identify high-order interactions in coupled dynamical systems. As an advancement, \cite{r40} further developed the adaptive signal lasso (ASigL) method, which has the desirable property of shrinking parameters fully to either 0 or 1, thus improving accuracy in network reconstruction.

 In this paper we propose a new version of signal lasso based on two kinds of penalty function to estimate the signal parameter and uncovering network topology in complex network with a small amount of observations. We find the tuning parameter can be set to a large enough values such that the signal parameter can be completely shrunk either 0 or 1, therefore there is no need to tuning parameter.  The theoretical properties of this method has been studied, and some simulations results are conducted to make fully comparison for all methods, which including in linear regression models with different assumption, evolutionary game dynamic model and synchronization model. Finally three dataset from the human experiment are deeply analyzed. All simulation results in terms of our proposed method is effective to uncover the signals presented in the models and greatly decreasing computational complexity of the method.

\section{Signal Lasso with additive penalty}\label{Sec:2}

Consider the general linear regression model
\begin{equation}\label{eq10}
{\bf Y}={\bf X} \beta+\epsilon
\end{equation}
where $\epsilon$ is a noise or random error with mean zero and finite variance, ${\bf X} = [x_{ij}]$ is an $n\times p$ matrix, ${\bf Y} = [y_i]$
is an $n\times 1$ vector, and $\beta = [\beta_i]$ is a $p \times 1$ unknown vector. $\epsilon$ is a random error.  To eliminate the intercept from (\ref{eq10}), throughout this paper, we centre the response variable and predictor variable so that the mean of response is zero. We assume the parameter $\beta$ here have a signal property, e.g. the true values of $\beta_j$ for $j=1,\cdots, p$ either 0 or 1. This kind of problem is very common in reconstruction of complex network to identifying the signal of connection or not connection (\cite{r14}; \cite{r26}).

In many practical situations, the signal always sparse, therefore lasso method, as a useful shrinkage estimation, have been used to find the estimator of parameter $\beta$ (\cite{r14}; \cite{r26}). 
Although  lasso method, have ability of shrink the parameter estimates toward to zero under the natural sparsity in complex network, the existent links between nodes can not be shrunk to its true value of 1, which will inevitably decrease the estimation accuracy in most cases. For this reason, \cite{r39} proposed signal lasso (SigL) method to solve the network reconstruction problem and they found SigL performed better than lasso and CS methods. However signal lasso still have some drawback in which the elements of $\hat\beta$ that fall in interval (0.1,0.9) can not be successfully selected to the correct class. Furthermore \cite{r40} proposed adaptive signal lasso, where appropriate weights are imposed to the penalty to promoting the estimation accuracy. The objective function with weighted penalty $\ell$ is defined by
\begin{equation}\label{eq12}
\ell(\beta,\lambda_1,\lambda_2\vert Y,X)=\frac{1}{2}||Y-X\beta||_2^2
+\lambda_1\sum_{j=1}^p\omega_{1j}\vert \beta_j\vert + \lambda_2\sum_{j=1}^p\omega_{2j}\vert \beta_j -1\vert,
\end{equation}
where $\omega_{1j}$ and $\omega_{2j}$ are two weights. The estimator defined by $\hat\beta^*=\arg\min_{\beta} \ell (\beta,\lambda_1,\lambda_2\vert Y,X)$ is called by adaptive signal Lasso(ASigL). The penalty function here is defined as 
$$PF(\beta)=\lambda_1 \sum_{j=1}^p \omega_{1j}\vert \beta_j\vert+\lambda_2 \sum_{j=1}^p \omega_{2j}\vert \beta_j-1\vert$$
They employed $\omega_{1k}=1$ and $\omega_{2k}=\vert\hat\beta_{k0}\vert$  for $k=1,\cdots,p$ in their study, where $\hat\beta_{k0}$ is an initial estimator of $\beta_{k}$, 
for example $\hat\beta_{k0}$ can be ordinary least square estimator for $p<n$ or ridge estimator for $p>n$. 

If Eq. (\ref{eq12}) is re-parametrized  by $\lambda=\lambda_2$ and $\alpha=\lambda_1/\lambda_2$, they found following fact when  $ \lambda \to +\infty$ for fixed $\alpha$
\begin{equation}\label{eq15}
\tilde \beta_k \to \left\{
\begin{array}{l}
 1, \ \ \ \ \ \ \ \hat \beta_{k0} > \alpha , \\
 0, \ \ \ \ \ \ \ \hat \beta_{k0} < \alpha, 
\end{array}
\right.
\end{equation}
since $\alpha_2 \to \alpha$ and $\alpha_1\to \alpha$ when $\lambda \to +\infty$. This result indicates that if $\lambda$ is selected large enough, the estimators from adaptive signal Lasso can be completely shrunk to either 0 or 1 and thus remove the unidentified set which will be presented in signal lasso method \cite{r39}. Signal lasso (SigL) corresponds the result using penalty function Eq. \ref{eq12} with two weight functions equal to 1, in which the property presented in Eq. \ref{eq15} is not true. 
The choices of tuning parameters $(\lambda_1,\lambda_2)$ in SigL can use well-known cross-validation (CV) technique. 
In adaptive signal lasso method, They specify a large value for $\lambda$ and only tuning the parameter $\alpha$, which only involved one parameter and will reduce the computation greatly. 
Furthermore, since $\alpha$ represents  the proportion of data compressed to 0 in the interval (0, 1), it should be less than 1 and greater than 0. 

\begin{figure*}
\begin{center}
\centering{\includegraphics[width=0.80\textwidth]{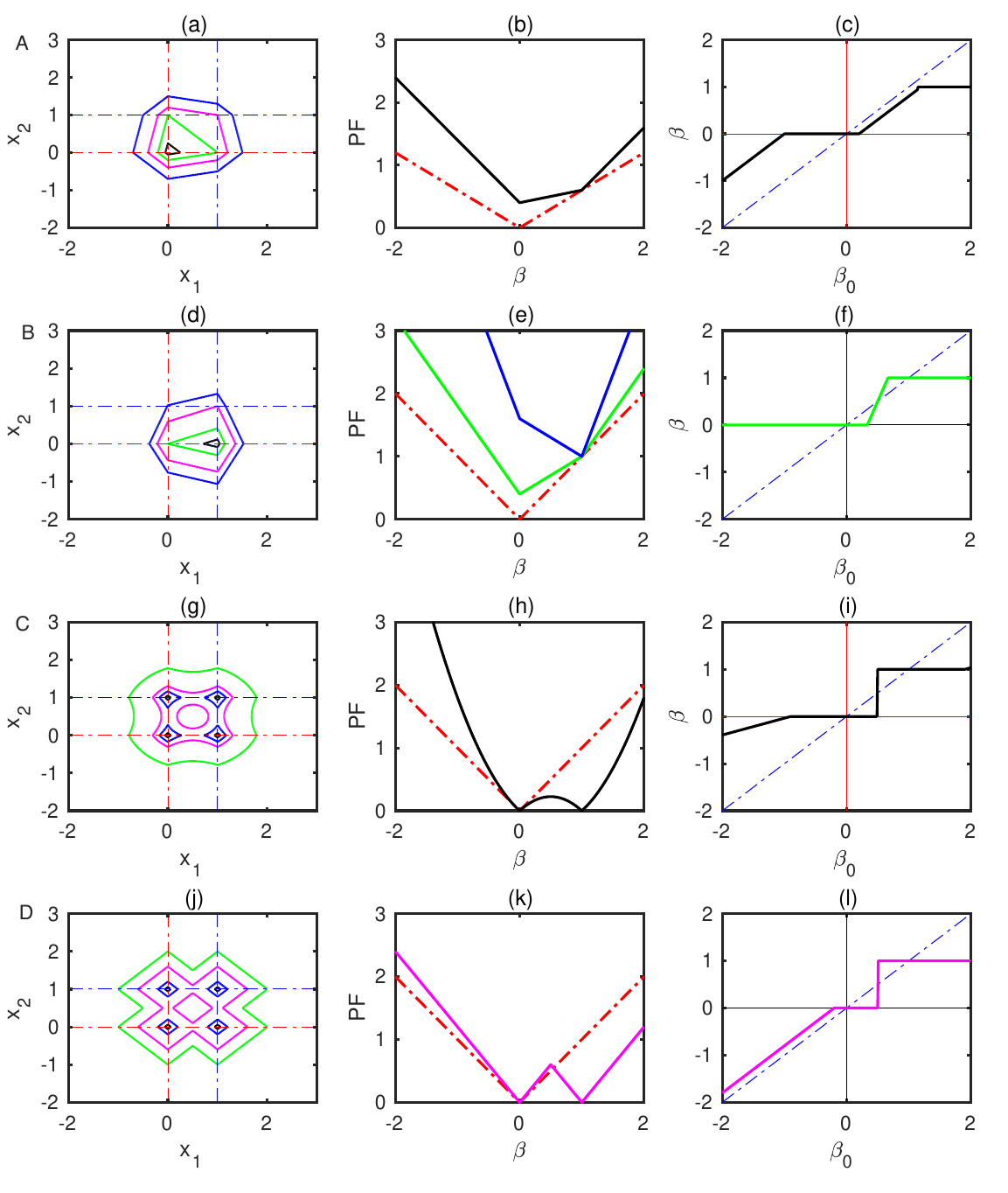}}
\end{center}
\caption{Some figures for four methods, where panel A, B, C and D are for signal lasso, adaptive signal lasso, signal lasso with product penalty, signal lasso with minimum penalty in two dimensional case (p=2), respectively.  
The first column represents the contour for penalty function of four methods, the second column represents the penalty function versus $\beta$, and last column represents the solution of $\beta$ under the orthogonal design. 
Constraint regions  are determined by $PF(x)=c$ for some constant $c=1,0.8,0.6,0.5$.  In second column, the red dashed line is the case of lasso method, and in Fig (e), $PF(x)= \vert x_1\vert+ \vert\hat\beta_{20} \vert \vert x_2-1 \vert $, where green line for case of $\hat \beta_{20}=0.2$ and blue line for case of $\hat \beta_{20}=0.8$. In Fig. (f) is the solution of $\beta$ under $\lambda_1=1,\lambda_2=2$. For signal lasso with production penalty, $\lambda=0.9$, while for signal lasso with minimum penalty, $\lambda=1.2$. }
\label{fig:1}
\end{figure*}

\section{Signal Lasso with non-convex penalty}\label{Sec:3}

\subsection{Method}
We first consider the linear regression mode (\ref{eq10}) with $cov(\epsilon)=\sigma^2 I_n$, the regression parameter $\beta_j$ have signal property that being either 0 or 1, therefore the accurate method to estimate $\beta$ is to minimize the following penalized least square

\begin{equation} \label{eq20}
L^*(\lambda, {\bf\beta})=\frac{1}{2} \sum_{i=1}^{n}(y_i-\sum_j {x}_{ij}{\beta}_j)^2+PF^*(\lambda,{\bf\beta}).
\end{equation}
Different from the previous penalties, here we consider two candidate of penalties
\begin{equation}\label{eq31}
	     PF^*(\lambda,\beta)=\lambda \sum_{i=1}^p \vert \beta_i (\beta_i-1) \vert 
\end{equation}
\begin{equation}\label{eq32}
	     	     PF^*(\lambda,\beta)=\lambda \sum_{i=1}^p min \{ \vert \beta_i \vert, \vert\beta_i-1\vert \} 
\end{equation}
Comparing with the additive penalty in signal lasso or adaptive signal lasso, these two penalties are not convex function, however they have appealing property that just one tuning parameter be involved and the $\beta_j$ can be shrunk to either 0 or 1 completely as shown later.   We will prove that there is no need to tuning parameter in our proposal, as long as  the regulation parameter selection is large enough, all the signals $\beta_j$ can be shrunk to either 0 or 1, which not only greatly reduce the computational complexity but also fully conforms to the  network reconstruction problem in Eq.\ref{eq20}.  

The product penalty function in Eq.\ref{eq31} has been used to identifying and estimating latent group structures in panel data (\cite{r41}), where they proposed classifier-Lasso (C-Lasso) to shrink individual coefficients to the unknown group-specific coefficients. For minimum penalty function in Eq.\ref{eq32}, \cite{r45} used similar functional form to study individualized model selection for different individuals and identify subgroups based on heterogeneous covariates' effects simultaneously. We use SL\_prod to denote signal lasso with product penalty and SL\_min for signal lasso with minimum penalty

Some patterns of penalty functions for SigL (panel A), ASigL (panel B), SL\_prod (panel C), SL\_min (panel D) in two dimensional case (p=2) are given in Fig. \ref{fig:1}, respectively. The first column represents the contour for penalty function of four methods, the second column represents the penalty function versus $\beta$, and last column represents the solution of $\beta$ under the orthogonal design.  Constraint regions  are determined by $PF(x)=c$ for some constant $c=1,0.8,0.6,0.5$. In second column, the red dashed line is the case of lasso method, and in Fig. \ref{fig:1} (e), $PF(x)= \vert x_1\vert+ \vert\hat\beta_{20} \vert \vert x_2-1 \vert $, where green line for case of $\hat \beta_{20}=0.2$ and blue line for case of $\hat \beta_{20}=0.8$. In Fig. \ref{fig:1} (f) is the solution of $\beta$ under $\lambda_1=1,\lambda_2=2$. For SL\_prod, $\lambda=0.9$, while for SL\_min, $\lambda=1.2$. 

It is clear that SigL method (panel A) fails to shrink the values of $\hat \beta$ that fall in the interval (0.1,0.9) to the correct class and leave an unclassified portion in network reconstruction. ASigL method improve this weakness such that the middle part between 0 and 1 can be shrunk toward to two directions. When one of tuning parameter is set to a large enough value, the slope in middle part of Fig. \ref{fig:1} (f) becomes a perpendicular line so that the parameters can be shrunk to either 0 or 1 (\cite{r40}). The signal lasso method with non-convex penalty in panel C and D of Fig. \ref{fig:1} obviously have this favorable property if tuning parameter $\lambda$ set a large value. The advantage of signal lasso with additive penalty (SigL and ASigL) is that their penalty functions are  convex and  the optimization problem does not suffer from  multiple local minima, and its global minimizer can be efficiently solved. However they need to select one or two tuning parameters which is computationally cost, especially for large size of network. Signal lasso method with non-convex penalty proposed in this paper do not need to tuning parameter as we show later, however the corresponding algorithm to solve the solution have to be developed.

\subsection{Updating formula for signal lasso with product penalty}
In this section, we propose a computationally very efficient algorithm for solving the two-directional signal lasso problem with product penalty (SL\_prod). Based on the coordinate descent method, the key step is to choose a single coordinate say $\beta_k$ at the $l$-th iteration and update $\beta_k^{(l+1)}$ by applying a univariate minimization over this coordinate, that is for $k=1, \ldots, p$
\begin{equation}\label{eq35}
\beta_k^{(l+1)}=\arg \min_{\beta_k}\frac{1}{2}\|\mathbf{Y}-\sum_{j< k}{\bf x}_j \beta_j^{(l+1)} -\sum_{j> k}{\bf x}_j \beta_j^{(l)}-{\bf x}_k \beta_k\|_2^2+\lambda |\beta_k||\beta_k-1|,
\end{equation}
where $\|\cdot \|_2^2$ is the $L^2$-norm of a vector. Let ${\bf\epsilon}_{k}^{(l)}=\mathbf{Y}-\sum_{j< k}{\bf x}_j \beta_j^{(l+1)} -\sum_{j> k}{\bf x}_j \beta_j^{(l)}$, then $r_k^{(l)}=\frac{{\bf x}_k^{\prime} {\bf\epsilon}_{k}^{(l)}}{{\bf x}_k^{\prime}{\bf x}_k} $ is the axis of symmetry of quadratic function part in Eq.(\ref{eq35})  that related to variable ${\bf x}_k$. Through some simple calculations, the update of  $\beta_k^{(l+1)}$ is obtained as
\begin{equation} \label{eq36}
\beta_k^{(l+1)}=\left\{ \begin{array}{ll}
 \min\left(0, \frac{{\bf x}_k^{\prime} {\bf\epsilon}_{k}^{(l)} +\lambda/2}{{\bf x}_k^{\prime}{\bf x}_k+\lambda} \right) & \text{if} \quad r_k^{(l)}\le 0 \\
 \max\left(0,  \frac{{\bf x}_k^{\prime} {\bf\epsilon}_{k}^{(l)} -\lambda/2}{{\bf x}_k^{\prime}{\bf x}_k-\lambda}  \right)\cdot \text{I}({\bf x}_k^{\prime}{\bf x}_k>\lambda) & \text{if} \quad 0<r_k^{(l)}\le 1/2 \\
  \min\left(1,  \frac{{\bf x}_k^{\prime} {\bf\epsilon}_{k}^{(l)} -\lambda/2}{{\bf x}_k^{\prime}{\bf x}_k-\lambda}  \right)\cdot \text{I}({\bf x}_k^{\prime}{\bf x}_k>\lambda)+\text{I}({\bf x}_k^{\prime}{\bf x}_k\le \lambda) & \text{if} \quad 1/2<r_k^{(l)}\le 1 \\
  \max\left(1, \frac{{\bf x}_k^{\prime} {\bf\epsilon}_{k}^{(l)} +\lambda/2}{{\bf x}_k^{\prime}{\bf x}_k+\lambda} \right) & \text{if} \quad r_k^{(l)}>1
\end{array}	\right. ,
\end{equation}
where $\text{I}(\cdot)$ is the indicative function. 

\subsection{Updating formula for signal lasso with minimum penalty}

The computational algorithm for solving the SL\_min can also be derived similarly. Based on the coordinate descent method, the key step is to choose a single coordinate say $X_k$ at the $l$-th iteration and update $\beta_k^{(l+1)}$ by applying a univariate minimization over this coordinate, that is for $k=1, \ldots, p$
\begin{equation}\label{eq38}
 \beta_k^{(l+1)}=\arg \min_{\beta_k}\frac{1}{2}\|\mathbf{Y}-\sum_{j< k}{\bf x}_j \beta_j^{(l+1)} -\sum_{j> k}{\bf x}_j \beta_j^{(l)}-{\bf x}_k \beta_k\|_2^2+\lambda \min\left(|\beta_k|, |\beta_k-1| \right),
\end{equation}
where $\|\cdot \|_2^2$ is the $L^2$-norm of a vector. Let ${\bf\epsilon}_{k}^{(l)}=\mathbf{Y}-\sum_{j< k}{\bf x}_j \beta_j^{(l+1)} -\sum_{j> k}{\bf x}_j \beta_j^{(l)}$, then $r_k^{(l)}=\frac{{\bf x}_k^{\prime} {\bf\epsilon}_{k}^{(l)}}{{\bf x}_k^{\prime}{\bf x}_k} $ is the axis of symmetry of quadratic function part in Eq.(\ref{eq38}) that related to variable ${\bf x}_k$. Through some simple calculations, the update of  $\beta_k^{(l+1)}$ is obtained as
\begin{equation} \label{eq39}
	\beta_k^{(l+1)}=\left\{ \begin{array}{ll}
		\text{sign}(r_k^{(l)})\cdot \max\left(0, |r_k^{(l)}|-\frac{\lambda}{{\bf x}_k^{\prime}{\bf x}_k} \right) & \text{if} \quad |r_k^{(l)}|\le|r_k^{(l)}-1|  \\
		1+\text{sign}(r_k^{(l)}-1)\cdot \max\left(0, |r_k^{(l)}-1|-\frac{\lambda}{{\bf x}_k^{\prime}{\bf x}_k} \right) & \text{if} \quad |r_k^{(l)}|>|r_k^{(l)}-1|
	\end{array}	\right. ,
\end{equation}
where $\text{sign}(\cdot)$ is the signum function. 

\subsection{Selection of tuning parameter}

{\bf Lemma}:  {\it For coordinate descent estimator of signal parameters with non-convex penalty given in (\ref{eq31}) or (\ref{eq32}), the estimator iteratively obtained from equation (\ref{eq36}) or (\ref{eq39}) can converge to either 0 or 1 when $\lambda\to\infty$.   }

{\it Proof: } The results can be proved from the formulae given in  equation  (\ref{eq36}) or (\ref{eq39}). For estimator using SL\_prod 
 when  $\lambda \to \infty$ , 
 $$ \frac{{\bf x}_k^{\prime} {\bf \epsilon}_{k}^{(l)} +\lambda/2}{{\bf x}_k^{\prime}{\bf x}_k+\lambda}\to 1/2,  \frac{{\bf x}_k^{\prime} {\bf \epsilon}_{k}^{(l)} -\lambda/2}{{\bf x}_k^{\prime}{\bf x}_k-\lambda}\to 1/2, \  \text{I}({\bf x}_k^{\prime}{\bf x}_k>\lambda)\to 0, $$ 
 thus $\hat\beta_k^{(l+1)}$ tends to 0 for $r_k^{(l)}\le 1/2$ and 1 otherwise. Finally the estimation of $\beta_k^{(l+1)}$ will be either 0 or 1.  
 For estimator using SL\_min, from equation (\ref{eq39}),  
 it is clear that $|r_k^{(l)}|-\frac{\lambda}{{\bf x}_k^{\prime}{\bf x}_k} $ and $|r_k^{(l)}-1|-\frac{\lambda}{{\bf x}_k^{\prime}{\bf x}_k}$ will eventually be less than 0, and finally the estimation of $\beta_k^{(l+1)}$ will be either 0 or 1. 
 
 This result indicates that we can set $\lambda$ by a large enough values and without tuning $\lambda$ using cross validation, and the signal parameter can be shrunk to 0 or 1 which achieve the purpose of complete classification. The overall iterative coordinate descent algorithm to find the solution can be obtained using these formulae.

\section{The metrics of reconstruction accuracy}\label{sec:3d}

To measure the accuracy of network reconstruction, we have to define some metrics for assessing the efficiency of proposed method. 
Follows the common strategy in network reconstruction in literature (\cite{r26}; \cite{r14}), the signal parameter $\beta$ can be classified as  signal ($\beta=1$) if $\hat \beta\in 1\pm 0.1$ and non-signal ($\beta=0$) if $\hat \beta\in 0\pm 0.1$, and remains are unclassified. Therefore there exist a portion of unclassified set if estimated $\hat\beta$ does not belong to either of two intervals. All quantities involved in classification can be written in Table \ref{tab: illus} (see \cite{r40}). 

When the predicted class can be completely classified to two classes, the most common used indexes for measuring accuracies include true positive rate (TPR, sensitivity or recall), true negative rate (TNR, or specificity) and precision (Positive prediction value, PPV ),  as well as AUROC (the area under the receiver operating characteristic curve) and AUPR (the area under the precision recall curve) (\cite{r20}; \cite{r14}; \cite{r46}), where TPR and TNR are defined by 
\begin{equation}\label{eq45}
TPR=\frac{TP}{TP+FN}, \ \ \ \  TNR=\frac{TN}{TN+FP},
\end{equation}
where $TP (TN)$ means the number that  signal (non-signal ) is correctly identified, $FP$ is the number that non-signal is falsely predicted as signal and  $FN$ is the number that signal is falsely predicted as non-signal. Among these measures,  Matthews correlation coefficient ($MCC$) is a measure that accounting for the effects of unbalanced number of signal and non-signal parameters (Chicco, 2017; Chicco and Jurman, 2020), and widely used in machine learning research.

To considered the effects of un-classification in Table \ref{tab: illus},  the success
rates for the detection of existing links (SREL) and non-existing links (SRNL) can be  defined to study the performance of network reconstruction (\cite{r26}; \cite{r14}; \cite{r39}). Another more efficient measure in this case to study the estimated accuracy when unclassified set remain is the adjusted $MCC (MCCa)$  defined by \cite{r40}, which has form of
\begin{equation}\label{eq47}
MCCa=\frac{TP\times TN-FPa\times FNa}{\sqrt{(TP+FPa)(TP+FNa)(TN+FPa)(TN+FNa)}}
\end{equation}
where $FNa = FN+UCP$ and $FPa = FP+UCN$.  Its values range from -1 to 1 and the large value indicates the good performance. It is clear when un-classified class disappears,$MCCa$ reduces to $MCC$. It is easy to see that MCCa play the similar role as $MCC$ when un-classification appear. \cite{r40} found $MCCa$ performs well in network reconstruction, thus we will use this metric to measure the accuracy of the method in this paper.

\begin{table*} \footnotesize
	\caption{Measures for accuracy of network reconstruction }
	\label{tab: illus}
	\begin{center}
		\begin{tabular*}{400pt}{@{\extracolsep\fill}cccccccccccc@{\extracolsep\fill}} \toprule
			\multirow{2}*{} &\multirow{2}*{Actual class }&\multicolumn{9}{@{}c@{}}{Predicted class } \\ 
			\cline{3-11} 
			&&\multicolumn{3}{@{}c@{}}{Signal }&\multicolumn{3}{@{}c@{}}{Non signal }&\multicolumn{3}{@{}c@{}}{Unclassified}\\  \hline
			&signal class &   \multicolumn{3}{@{}c@{}}{True positive (TP) }&\multicolumn{3}{@{}c@{}}{False negative (FN) }&\multicolumn{3}{@{}c@{}}{Unclassified positive (UCP)} \\
			&Non signal class &   \multicolumn{3}{@{}c@{}}{False positive (FP) }&\multicolumn{3}{@{}c@{}}{True negative (TN) }&\multicolumn{3}{@{}c@{}}{Unclassified negative (UCN)} \\
			\hline
		\end{tabular*}
	\end{center}
\end{table*}

\section{Numerical Studies}

In this section, we use two kinds of model to conduct the simulation studies. One based on standard linear regression model with different assumptions, which including different error distribution and correlated or uncorrelated design matrix. Second one use the network reconstruction model based on game-based evolutionary dynamics and Kuramoto model as used in \cite{r39, r40}.

\subsection{Linear regression models}

The model generation is given by
\begin{equation}\label{eq50}
Y=\bf {1}_n \beta_0+X_1 \beta_1+X_2\beta_2+\epsilon
\end{equation}
where $\beta_0=0$ denote the intercept, $\bf{1}_n$ is a $n\times 1$ vector with all elements equal to 1, $\beta_1 \in R^{p_1}$ denote the signal parameter with elements 1, and $\beta_2 \in R^{p_2}$ denote the non-signal parameter with elements 0, and  $\epsilon$ is a noise term. Smaller $p_1$ is called by sparse signal and larger $p_1$ (comparing with n and $p=p_1+p_2$) is called dense signal. We consider several case in simulations.

\begin{itemize}
  \item \textit{Case I}: Each design matrix come from standard normal score with mean zeros and variance 1, but columns in in $X$ are dependent with the correlation between $x_i$ and $x_j$ is $r^{|i-j|}$ with $r=0.5$ (see \cite{r25} and \cite{r1} for similar settings). The error variable $\epsilon$ is generated from normal distribution with mean zeros and variance $\sigma=0.4, 1, 2$ respectively. Two combination of sample size and number of independent variable are considered by $n=100,p=30$ and $n=50, p=150$, respectively. The sparse signal with $p1=6$, non-sparse signal with $p1=20$ and dense signal with $p1=100$ are considered for $n<p$.
  \item \textit{Case II}: The design matrix are the same as Case I, but the error variable are generated from exponential distribution $Exp(\sigma)$, and Gamma distribution $Ga(4,\sigma/2)$, $(n, p)$ are the same as \textit{Case I}. The exponential and gamma distributions are independently generated. Other settings are same as Case I. 
  \item \textit{Case III}: The design matrix are the same as Case I, but the error variable $\epsilon$ has a AR(1) with correlation $\rho$ ($\epsilon_i=\rho \epsilon_{i-1}+u_i$, where $u_i$ is a Gaussian random variable with mean zero and variance $\sigma^2$ ), and  $\rho$ takes values of 0.4 or 0.7 respectively. 
\end{itemize}

These three cases corresponds to different sources of the data come from. Case I is for situation of standard linear model with error variable has a independent normal distribution, case II is for case when data has non-normal distribution and case III represents the scenario of correlated observations. After getting the simulated data set, the signal parameter $\beta$ can be classified as  $\beta=1$ if $\hat \beta\in 1\pm 0.1$ and $\beta=0$ if $\hat \beta\in 0\pm 0.1$, and remains are unclassified (\cite{r26}; \cite{r14}).  Some indices in classification problem, including TPR, TNR and Precision rate(PCR), SREL and SRNL can be calculated to quantitatively measure the performance (\cite{r20}; \cite{r14}; \cite{r39, r40}). However for saving the space of the paper, we only list several important measures such as MSE, UCR, MCC and MCCa in our simulations. 

\begin{table*}[htbp]\footnotesize
	\begin{center}
		\caption{Simulation results in Case I for nine methods: Lasso, Adaptive Lasso, SCAD, MCP, Elastic Net, Signal Lasso, Adaptive Signal Lasso, Signal lasso with product penalty and minimum penalty in linear regression models. All of the results are averaged over 500 independent realizations, where n is the sample size, $p$ is the number of explanatory variables, $p_1$ is the number of signals (number of $\beta=1$). The first panel is for case of $p<n$ and signal is sparse with $ p_1=6$ . The second panel is for $p>n$  and sparse signal. The third panel consider the cases of non-sparse signal with $p_1=20$. The bottom panel consider the cases of dense signal with $p_1=100$. The noise is introduced by $\sigma = 0.4, 1, 2$, respectively, in three columns. }
		\label{tb:2} 
\setlength{\tabcolsep}{3mm}{
	\begin{tabular}{ccccccccc}	
			\hline
Method	&  \multicolumn{3}{c}{MSE/UCR/MCC/MCCa }  \\
\hline
$(n, p, p_1, \sigma)$ 		&$(100, 30, 6, 0.4)$ & $(100, 30, 6, 1)$  & $(100, 30, 6, 2)$\\	
\hline
Lasso &0.0014/0.027/1.000/0.919 & 0.0083/0.221/0.980/0.347& 0.0320/0.359/0.863/0.022 \\
A-lasso &0.0015/0.343/1.000/0.900 & 0.0095/0.232/0.985/0.316 &0.0459/0.349/0.757/-0.006\\
SCAD & 0.0008/0.023/1.000/0.937 &  0.0043/0.108/1.000/0.632&0.0509/0.226/0.702/0.150\\
MCP & 0.0008/0.019/1.000/0.943 &  0.0046/0.111/1.000/0.633 &0.0506/0.220/0.736/0.184\\
ElasticNet&0.0012/0.333/1.000/0.892 & 0.0066/0.139/0.970/0.503 &0.0239/0.201/0.721/0.225\\ 
SigL&0.0007/0.012/1.000/0.966& 0.0046/0.133/1.000/0.609& 0.0232/0.256/0.943/0.249\\
ASigL &0.0001/0.002/1.000/0.993 & 0.0006/0.012/1.000/0.965 &0.0086/0.016/0.982/0.933\\
SL\_prod &   0.0000/0.000/0.999/1.000 & 0.0000/0.000/0.999/1.000 &0.0058/0.000/0.982/0.982\\
SL\_min &   0.0000/0.000/0.999/1.000 & 0.0000/0.000/0.999/1.000 &0.0058/0.000/0.982/0.982\\    

\hline
$(n, p, p_1, \sigma)$ 	&	$(50, 150, 6, 0.4)$ & $(50, 150, 6, 1)$ & $(150, 150, 6, 2)$ \\
\hline
Lasso &0.0016/0.043/0.990/0.529 & 0.0085/0.162/0.970/0.107&0.034/0.268/0.559/-0.042\\
A-lasso &0.0026/0.076/0.980/0.375 & 0.0130/0.191/0.860/0.061&0.0415/0.264/0.577/0.035\\
SCAD & 0.0002/0.007/1.000/0.896 & 0.0111/0.037/0.609/0.305 &0.0315/0.043/0.141/0.038\\
MCP &0.0002/0.007/1.000/0.892 & 0.0130/0.032/0.588/0.309&0.0365/0.031/0.144/0.065\\
ElasticNet& 0.0007/0.019/1.000/0.712 & 0.0037/0.039/0.850/0.359& 0.0126/0.050/0.439/0.144\\
SigL &  0.0004/0.013/1.000/0.839 & 0.0028/0.051/1.000/0.482&0.0127/0.085/0.883/0.255 \\
ASigL &0.0009/0.001/0.989/0.976 & 0.0213/0.005/0.821/0.781&0.0111/0.003/0.873/0.841\\
SL\_prod &   0.0000/0.000/0.999/1.000 & 0.0030/0.000/0.967/0.967&0.0584/0.000/0.559/0.559\\
SL\_min &   0.0000/0.000/0.999/1.000 & 0.0030/0.000/0.967/0.967&0.0584/0.000/0.559/0.559 \\    

\hline

$(n, p, p_1, \sigma)$ 	&	$(50, 150, 20, 0.4)$ & $(50, 150, 20, 1)$ &	$(50, 150, 20, 2)$\\
\hline
Lasso &0.0279/0.202/0.905/0.091 & 0.0431/0.255/0.847/-0.01&0.0796/0.311/0.600/-0.108 \\
A-lasso &0.0104/0.140/0.975/0.372 & 0.0305/0.237/0.918/0.057&0.0802/0.285/0.646/-0.066 \\
SCAD & 0.1556/0.085/0.095/-0.01 &  0.1565/0.082/0.124/0.057&   0.1680/0.077/0.107 /0.004\\
MCP &0.1709/0.072/0.078/-0.01  & 0.1771/0.069/0.091/-0.00 &0.1919/0.058/0.064/-0.008\\
ElasticNet& 0.0142/0.119/0.961/0.362 & 0.0276/0.153/0.883/0.172 &0.0533/0.157/0.609/0.068\\
SigL &  0.0015/0.044/1.000/0.815 &0.0092/0.118/0.997/0.549 & 0.0349/0.161/0.927/0.370\\
ASigL &0.0368/0.005/0.840/0.819 & 0.0400/0.004/0.827/0.809 &0.0532/0.004/0.767/0.751\\
SL\_prod &   0.0028/0.000/0.988/0.988 & 0.0168/0.000/0.931/0.931 &0.0835/0.000/ 0.689/0.689\\
SL\_min &  0.0028/0.000/0.988/0.988 & 0.0168/0.000/0.931/0.931 &0.0835/0.000/0.689/0.689 \\   
   
\hline

$(n, p, p_1, \sigma)$ 	&	$(50, 150, 100, 0.4)$ & $(50, 150, 100, 1)$ &	$(50, 150, 100, 2)$\\
\hline
Lasso &0.7554/0.408/0.120/-0.382 & 0.7435/0.408/0.121/-0.376&0.7622/0.411/0.135/-0.376 \\
A-lasso &0.7498/0.364/0.122/-0.268 & 0.7388/0.361/0.135/-0.257&0.7671/0.361/0.146/-0.257 \\
SCAD & 0.9004/0.082/0.045/-0.058 &  0.8662/0.076/0.042/0.050&   0.8545/0.080/0.039 /-0.053\\
MCP &1.0681/0.046/0.015/-0.048  & 1.0106/0.045/0.007/-0.044 &0.9901/0.037/0.010/-0.045\\
ElasticNet& 0.5375/0.311/0.016/-0.185 & 0.5337/0.313/0.153/0.188 &0.5314/0.321/0.177/-0.175\\
SigL &  0.0208/0.145/0.981/0.657 &0.0414/0.202/0.966/0.527 & 0.0837/0.213/0.886/0.448\\
ASigL &0.2982/0.008/0.437/0.427 & 0.3080/0.007/0.425/0.416 &0.3049/0.006/0.432/0.424\\
SL\_prod &   0.1985/0.000/0.597/0.597 & 0.2082/0.000/0.582/0.582 &0.1968/0.000/ 0.601/0.601\\
SL\_min &  0.1985/0.000/0.597/0.597 & 0.2082/0.000/0.582/0.582 &0.1968/0.000/0.601/0.601 \\   
   
\hline
		\end{tabular}
		}
	\end{center}
\end{table*}

We first consider the Case I and the results are listed in Table \ref{tb:2}, where upper panel is for sparse signal with $n>p$, middle panel is for sparse signal with $n<p$, while bottom panel is about non-sparse signal with $n<p$. The three column in each pane give the results of different variance with $\sigma=0.4,1,2$, respectively. All measures as listed in \cite{r40} are calculated but here we only list the results of MSE, UCR, MCC and MCCa as we see these four measures can reveal the overall accuracy of network reconstruction (refer to \cite{r40}).  
In Table \ref{tb:2}, we calculate the results based on nine methods in which the firs five methods are well known as shrinkage estimators for comparisons, and the last four methods are called as signal-lasso-type estimator because they shrink the parameters towards either 0 or 1. The lasso estimator is used as the initial estimator in last four methods.  It is clear all results show that signal-lasso-type estimators superior to the first five well known methods in estimating signal parameter. Therefore we focus on comparison among the last four methods, and have following conclusion from Table \ref{tb:2}:

(1) In cases of $n>p$ and the signal is sparse (first panel), we find that the SL\_prod and SL\_m outperform SigL and ASigL for different variance. In such cases, the SigL is not performs well, especially in case of large variance. When signal is non-sparse or dense, the results are the same (no show here).

(2) For the scenario of $n<p$ and sparse signal ($p_1=6$) (second panel), when variance of error variable is not large such as $\sigma =0.4$ and $1$, SL\_prod and SL\_m
clearly outperforms the ASigL and SigL based four measures. However when variance  is large such as $\sigma=2$, ASigL performs best with large MCCa and small MSE, although its UCR is not exactly equal to 0. In this case, the SigL performs poorly, especially for large variance. 

(3) For the scenario of $n<p$ and non-sparse signal ($p_1=20$) (third panel), SL\_prod and SL\_min performs best and followed by ASigL method. However for large variance with $\sigma=2$, it is observed that ASigL outperforms other three methods, where SigL is not  good. The results from second and third panel indicate that AsigL is robust against  variance of distribution, which again evidence the conclusion of \cite{r40}.

(4) Last panel of Table \ref{tb:2} list the results for case of  $n<p$ and dense signal ($p_1=100$). It is of interest to find that SigL become better for small variance with smallest MSE and largest MCCa and MCC even its UCR is not closed to zero. When variance tends to large, the advantage of SL\_prod and SL\_min standout, where they have largest MCCa and small UCR. In this case the SigL have smallest MSE and largest MCC but non-zero UCR. Therefor the performance of SigL is comparable with SL\_prod and SL\_min for  case of dense signal. The ASigL dose not perform well in the circumstance of dense signal.

\begin{figure*}
\centering{\includegraphics[width=0.8\textwidth]{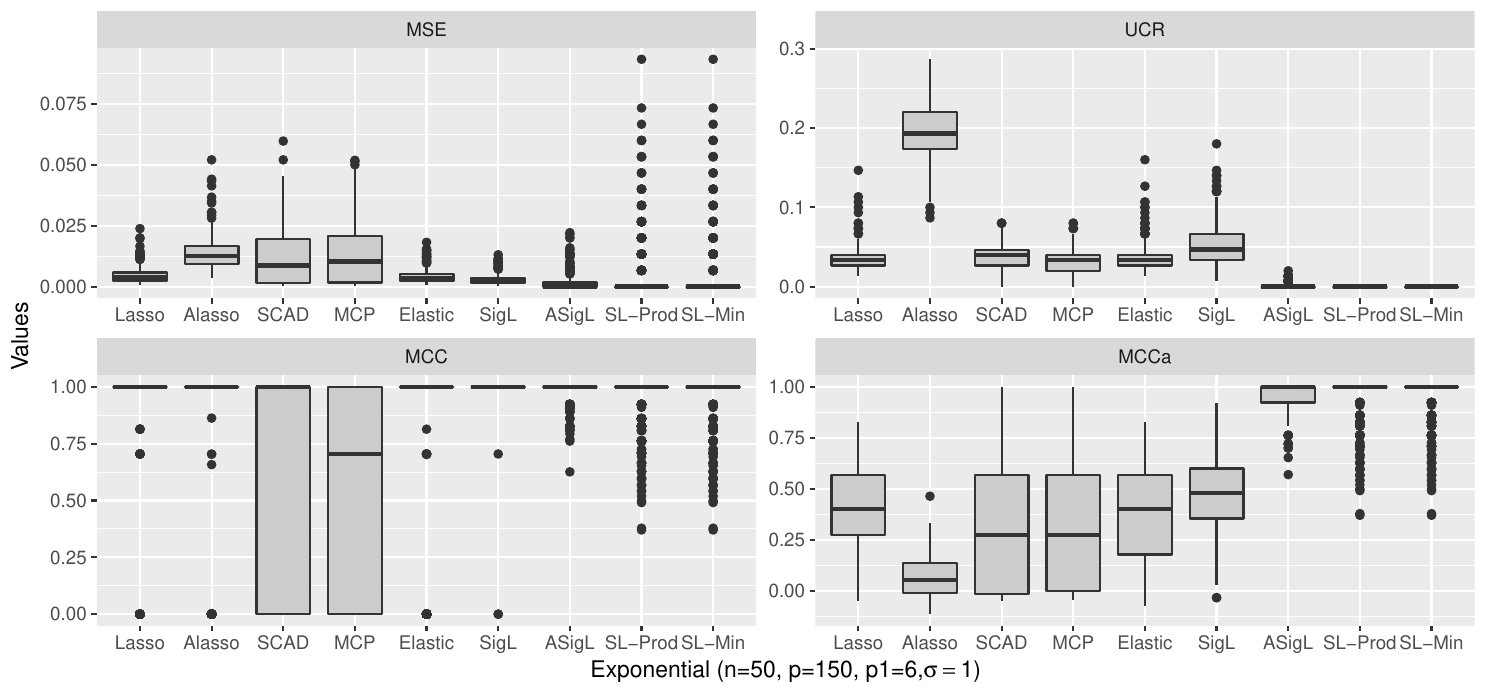}}
\caption{Boxplots of MSE, UCR, MCC and MCCa in Case II under exponential distribution for nine methods: Lasso, Adaptive Lasso, SCAD, MCP, Elastic Net, SigL, ASigL, SL\_prod and SL\_min in linear regression models. All of the results are averaged over 500 independent realizations, where n is the sample size, $p$ is the number of explanatory variables, $p_1$ is the number of signals (number of $\beta=1$). We consider case of $p>n$ and signal is sparse with $ p_1=6$. The noise distribute as exponential with $\sigma = 1$. }
\label{fig:2}    
\end{figure*}

\begin{figure*}
\centering{\includegraphics[width=0.8\textwidth]{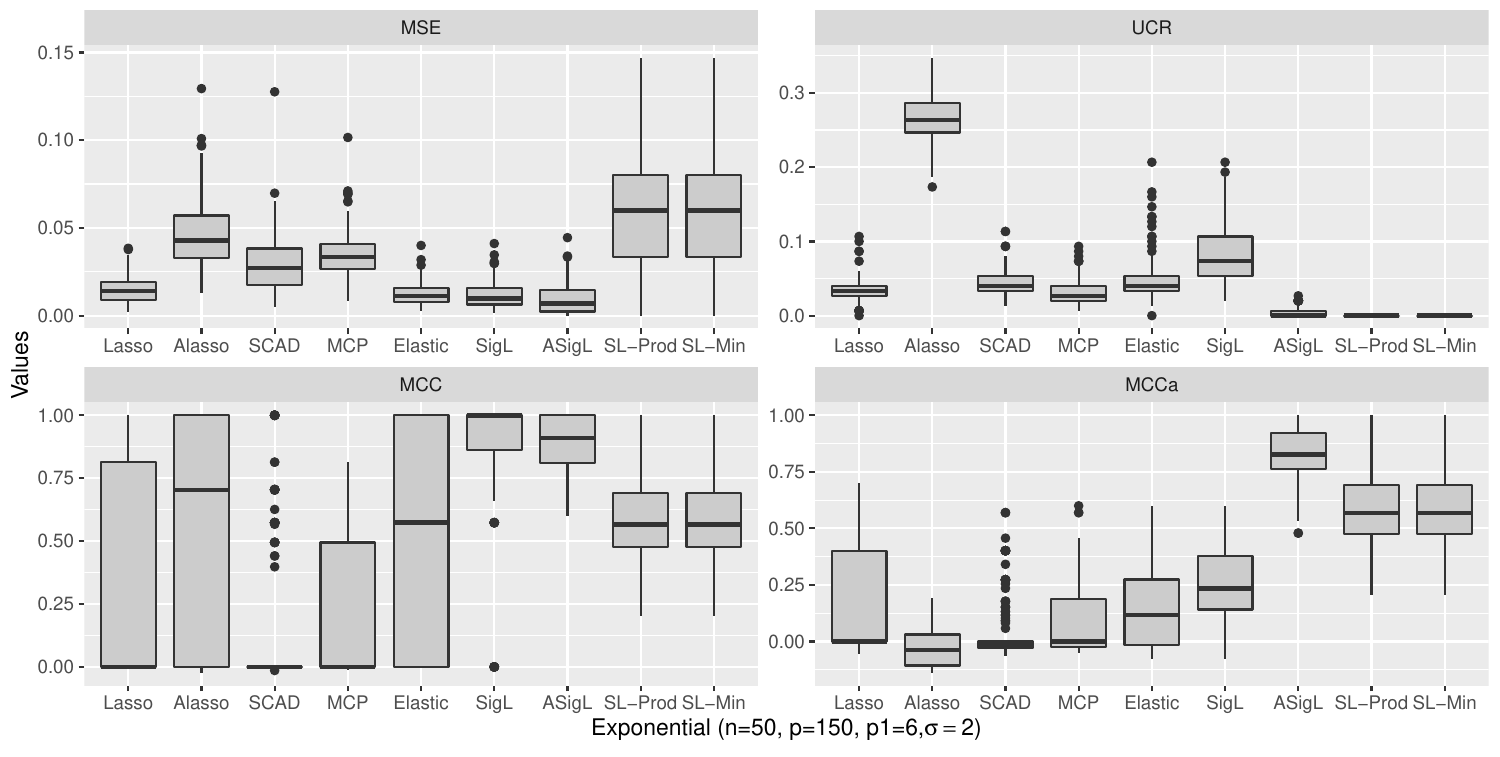}}
\caption{Boxplots of MSE, UCR, MCC and MCCa in Case II under Gamma distribution for nine methods: Lasso, Adaptive Lasso, SCAD, MCP, Elastic Net, SigL, ASigL, SL\_prod and SL\_min in linear regression models. All of the results are averaged over 500 independent realizations, where n is the sample size, $p$ is the number of explanatory variables, $p_1$ is the number of signals (number of $\beta=1$). We consider case of $p>n$ and signal is dense with $ p_1=20$. The noise distribute as exponential with $\sigma = 2$. }
\label{fig:3}    
\end{figure*}

Next we look at the Case II with Exponential or Gamma distribution are assumed to error variable.  We find the results are similar to those of Gussian distribution. We only list the results for exponential distribution. Figure (\ref{fig:2}) is the plots for MSE, UCR, MCC and MCCa for nine methods, shown that SL\_prod and SL\_min have larger MCCa and MCC and smaller MSE and UCR than ASigL. However if variance $\sigma$ equal to 2, ASigL is better comparing with SL\_prod and SL\_min as shown in  Figure (\ref{fig:3}). The box-plot indicates that the SL\_prod and SL\_min have more outliers in MCCa, MCC and MSE than that of ASigL. The results based on Gamma distribution have the very similar patterns and the results are omitted here. 

\begin{figure*}
\centering{\includegraphics[width=0.8\textwidth]{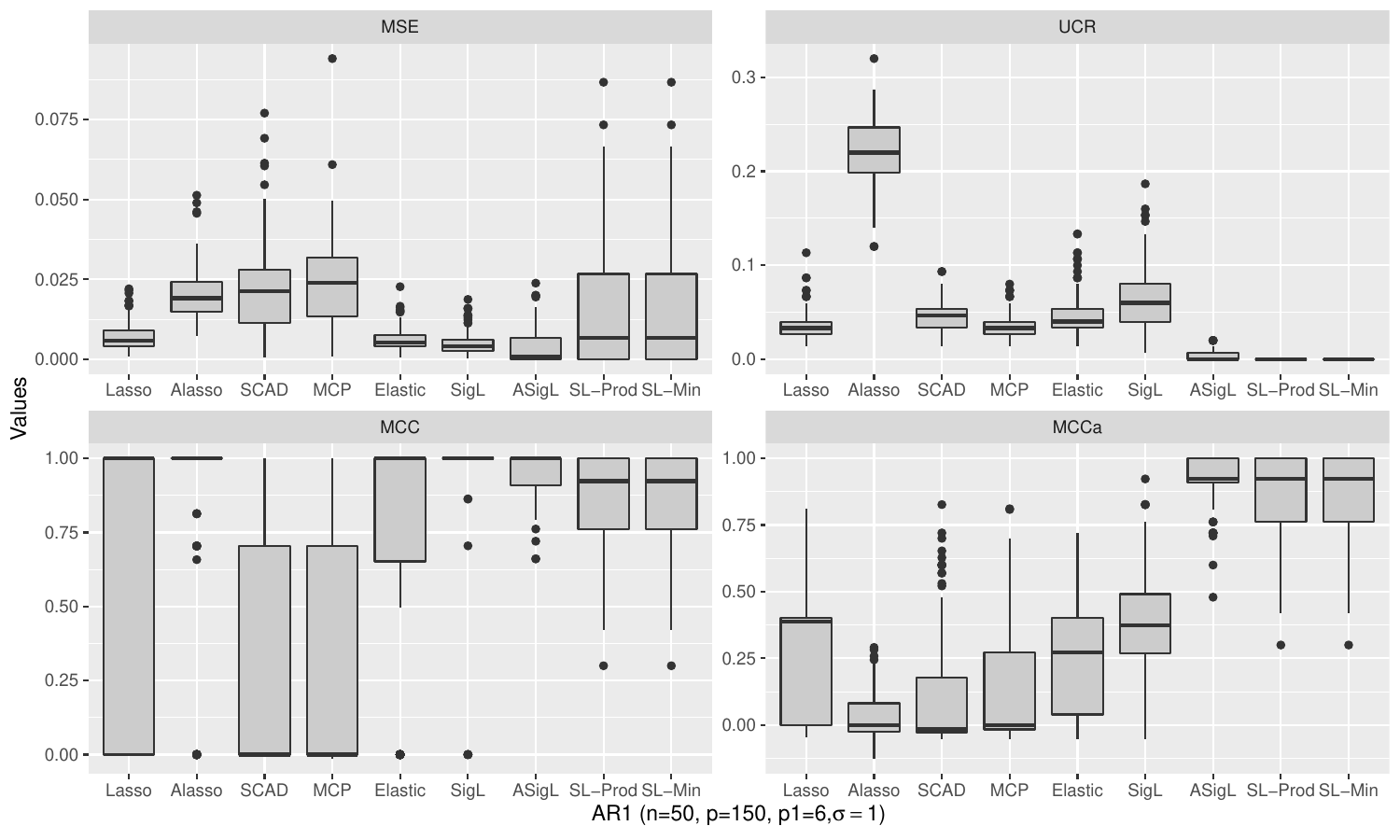}}
\caption{Boxplots of MSE, UCR, MCC and MCCa in Case III under AR(1) distribution for nine methods: Lasso, Adaptive Lasso, SCAD, MCP, Elastic Net, SigL, ASigL, SL\_prod and SL\_min in linear regression models. All of the results are averaged over 500 independent realizations, where n is the sample size, $p$ is the number of explanatory variables, $p_1$ is the number of signals (number of $\beta=1$). We consider case of $p>>n$ and signal is sparse with $ p_1=6$. The noise distribute as AR(1) with $\sigma = 1$. }
\label{fig:4_1}    
\end{figure*}

\begin{figure*}
\centering{\includegraphics[width=0.8\textwidth]{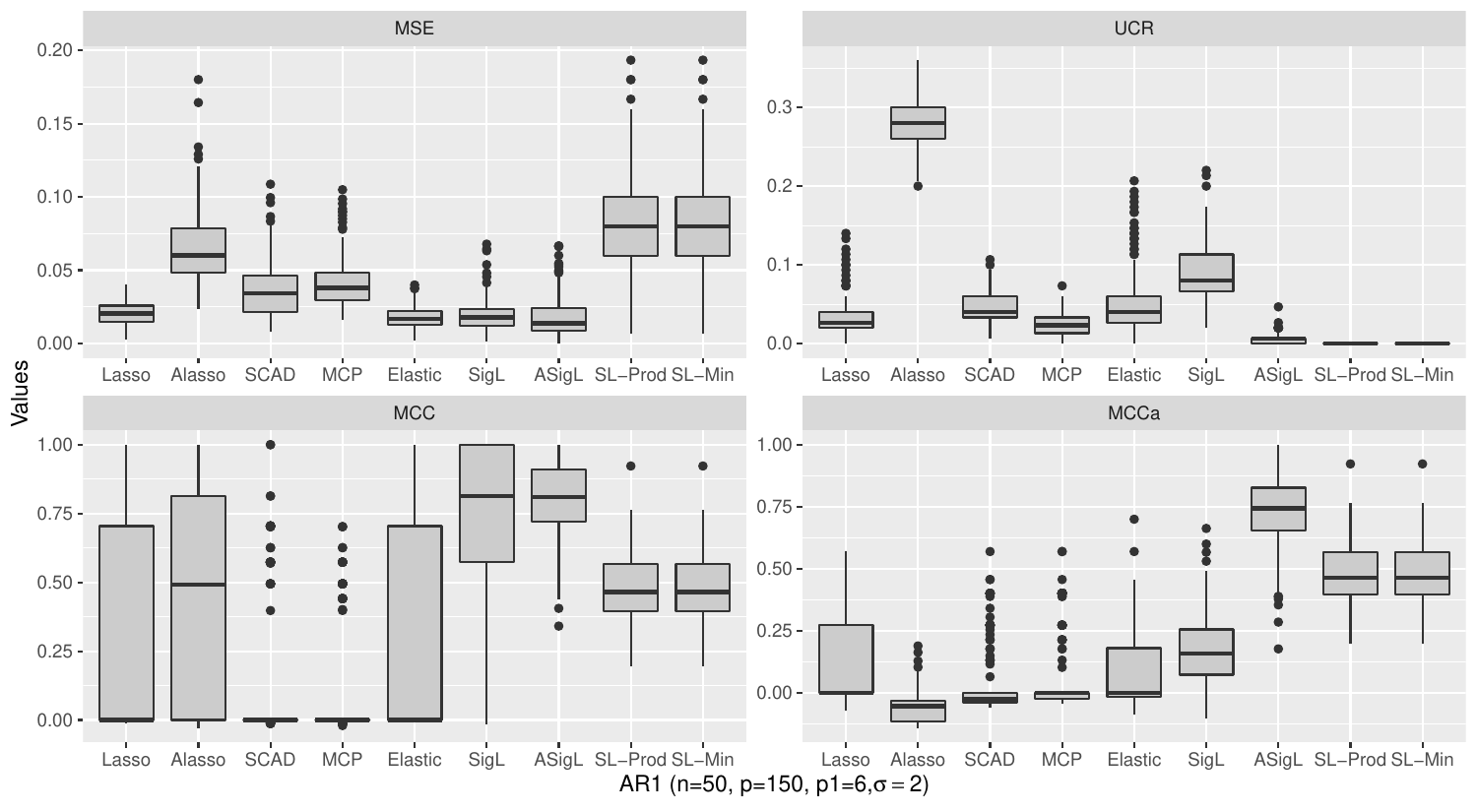}}
\caption{Boxplots of MSE, UCR, MCC and MCCa in Case III under AR(1) distribution for nine methods: Lasso, Adaptive Lasso, SCAD, MCP, Elastic Net, SigL, ASigL, SL\_prod and SL\_min in linear regression models. All of the results are averaged over 500 independent realizations, where $n$ is the sample size, $p$ is the number of explanatory variables, $p_1$ is the number of signals (number of $\beta=1$). We consider case of $p>>n$ and signal is sparse with $ p_1=6$. The noise distribute as AR(1) with $\sigma = 2$. }
\label{fig:4_2}    
\end{figure*}

Finally we consider Case III with linear regression with $n=50, p=150, p_1=6$, but error variable have correlated structure of AR(1) distribution with $\rho=0.7$.  When variance is smaller with $\sigma=1$, Figure (\ref{fig:4_1}) show that SL\_prod and SL\_min performs better than ASigL with larger MCCa and smaller MSE and UCR. However when $\sigma=2$, Figure (\ref{fig:4_2}) clearly indicates that ASigL outperforms SL\_prod and SL\_m. For non-sparse signal and dense signal, the results are very similar to those based on Gaussian distribution and the results are not shown here. 

In summary, we find that signal-lasso-type method performs well than the first five methods, where signal lasso with non-convex penalties seems overall outperform other methods in most cases. In the scenario of sparse and non-sparse signal, for small variance of noise, SL\_prod and SL\_min perform better than adaptive signal lasso, and they have exact zero UCR for Gaussian and other distribution, even for AR(1) correlated structure. However for large noise variance, ASigL outperform to SL\_prod and SL\_m, shown a better robustness for the noisy data. However for dense signal, our finding reveal that SigL performs better for small variance, while SL\_prod and SL\_min performs better for large variance. The results also show that the signal-lasso-type methods are robust against distribution of the data. 

\subsection{\label{sec3}Evolutionary-game-based dynamical model}

In an evolutionary game on structured populations, there are multiple players participate the game with their neighbors, where a node represents a player, and a link indicates that two players have a game relationship. The prisoner's dilemma game (PDG), snowdrift game (SDG), or spatial ultimatum game (SUG) can be used for network reconstruction (\cite{r26}; \cite{r14}; \cite{r39}). 
We illustrate our method by iterative game dynamics through Monte Carlo simulation.  The payoff matrix is given in Eq. (\ref{eq1}), here we take a simple structure with $R=1,T=b=1.15$, and $P=S=0$ in our simulation. 
The game  iterates forward in a Monte Carlo manner and player $i$ (the focal player) acquires its fitness (total payoff) $F_i$ by playing the game with all its direct neighbors, i.e., $F_i=\sum_{j=1}^N a_{ij} P_{ij}$. The focal player then randomly picks a neighbor $j$, which similarly acquires its fitness. For simulation cases, player $i$ tries to imitate the strategy of player $j$ with  Fermi updating probability $W=1/(1+\exp [(F_i-F_j)/K])$, where $K=0.1$ (\cite{r42}; \cite{r43}). By this updating rule of strategies, a series dataset can be generated for each player. To make the model more realistic, we account for a mutation of small rate for the data.  

Now $F_i=\sum_{j=1, j \neq i}^N a_{ij} P_{ij}$ can be written as a linear regression model 
\begin{equation}\label{eq65}
Y_i={\bf X}_i \tilde \beta_i +\epsilon_i,
\end{equation}
where $Y_i=(F_i(t_1), F_i(t_2), \cdots, F_i(t_L))' $, $\tilde \beta_i=(a_{i1}, \cdots, a_{iN})'$, and $/X_i$ has the form of 
\begin{displaymath}
\left(\begin{array}{cccccc}
P_{i1}(t_1) & \cdots & P_{i,i-1}(t_1) & P_{i,i+1}(t_1) & \cdots & P_{i_N}(t_1)\\
P_{i1}(t_2) & \cdots & P_{i,i-1}(t_2) & P_{i,i+1}(t_2) & \cdots & P_{i_N}(t_2)\\
\vdots & \vdots & \vdots & \vdots      \\
P_{i1}(t_L) & \cdots & P_{i,i-1}(t_L) & P_{i,i+1}(t_L) & \cdots & P_{i_N}(t_L)
\end{array}\right), 
\end{displaymath}
Let $Y=(Y'_1, \cdots, Y'_N)'$, $\beta=(\tilde \beta'_1, \cdots, \tilde \beta'_N)'$, ${\bf X}=diag({\bf X}_1, {\bf X}_2, \cdots, {\bf X}_N)$, then  Eq. (\ref{eq65}) can be converted into the general form of Eq. (\ref{eq10}). 

Here, $L$ represents the length of the data and indicates that there are $L$ interactions between players in the game dynamics. In many real-world scenarios, one needs to uncover a graph's topology using a limited amount of data. The amount of data is typically scaled by $\Delta = L/N$.
In Eq. (\ref{eq65}), $\beta$ contains the elements of the connectivity matrix $A = [a_{ij}]$, representing a signal parameter where a value of 1 indicates a connection between two nodes, and 0 indicates no connection. To evaluate the accuracy of the estimator for $\beta$, we plot various metrics as a function of $\Delta = L/N$, which helps to determine how much data is needed to achieve high estimation accuracy. This approach has been widely used in network reconstruction (\cite{r26}; \cite{r14}; \cite{r39, r40}). 

\begin{figure*} 
\centering{\includegraphics[width=0.8\textwidth]{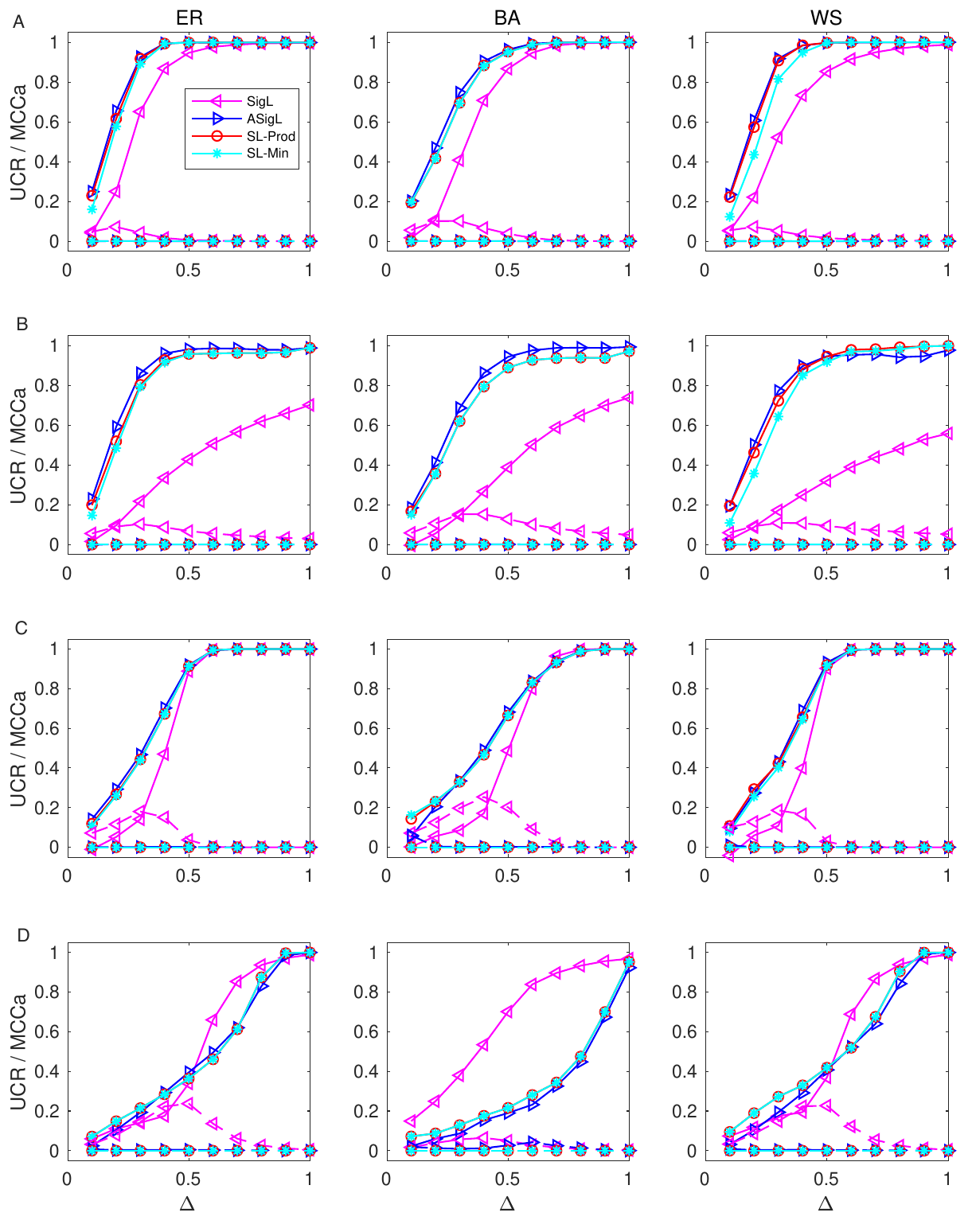}}
\caption{Accuracy measures MCCa and UCR in the reconstruction vs. $\Delta=L/N$, for PDG model attained by  method of SigL, ASigL, SL\_prod and SL\_min in three kinds of network. The panel A refers to results for PDG game with noise $\sigma^2=0.1$ and average degree 6. The panel B refers to the results for PDG game with noise variance $\sigma^2=0.3$ and average degree 6. The panel C refers to the results for PDG game with average degree 20 without noise. The panel D refers to the results of dense signal for PDG game with average degree 70 with variance of noise $\sigma^2=0.1$. Three columns correspond to the results based on ER random network,  BA scale-free network and  WS small world network, respectively. The network size $N=100$ and each point is averaged over 20 simulations. We only provide annotations for the graphical elements in the first figure in Panel A, and the others are the same, where the upper three curves are for MCCA and lower curves for UCR.  }
\label{fig:5}    
\end{figure*} 

Fig. \ref{fig:5} list the plot of accuracy measures against the scaled data length $\Delta=L/N$ in the PDG model with three networks under different conditions, using the methods of signal lasso, adaptive signal lasso and signal lasso with non-convex penalty. We only list these four methods for comparisons and other methods in Table \ref{tb:2} are ignored since SigL and ASigL outperform other methods as shown in \cite{r40}. We consider four cases with different assumptions for noise variances and connectivity density as shown in four panels, and three columns correspond to the results based on Erd\"os-R\'enyi (ER) random network,  Barab\'asi-Albert (BA) scale-free network and  small world (WS) network, respectively.

Panel A in Fig. \ref{fig:5} represents the case with $N=100$, an average degree of $k=6$, and a noise variance of $\sigma^2=0.1$. It is evident that SigL performs the worst across the three networks, showing lower reconstruction accuracy (MCCa) and higher UCR at small $\Delta$. Both ASigL and SL\_prod perform well, with ASigL slightly outperforming SL\_prod at small $\Delta = L/N$. However, this differs from the results in Table \ref{tb:2}, where the MCCa values obtained from SL\_prod are higher than those of SL\_min, particularly in the WS network. Upon examination, we find that this discrepancy arises because several columns in the design matrix $X$ are zero-columns (i.e., all elements in those columns are zero) due to the network's sparse connectivity. In such cases, the SL\_prod method is more robust than SL\_min, as demonstrated in the linear regression model in Table \ref{tb:A1} of Appendix. Panel B in Fig. \ref{fig:5} uses the same parameters as Panel A but with a higher noise variance of $\sigma^2=0.3$. The results are similar to those in Panel A in most cases, where SigL performs worse and ASigL outperforms the other three methods. However, in the WS network, SL\_prod and SL\_min have higher MCCa values than ASigL when $\Delta > 0.6$. Additionally, we also observe that SL\_prod outperforms SL\_min at small $\Delta = L/N$.

Panel C in Fig. \ref{fig:5} presents the results for the case with $N=100$ and an average degree of $k=20$, representing a non-sparse network. The performances of ASigL, SL\_prod, and SL\_min are similar, while SigL still performs poorly.
Panel D shows the results for a dense signal scenario with an average degree of 70 and a noise variance of $\sigma^2=0.1$. In this case, SigL outperforms the other three methods in the BA network. In both the ER and WS networks, the performance of SigL is comparable to the other three methods, showing smaller MCCa when $\Delta < 0.5$ but larger MCCa when $\Delta > 0.5$. However, it is notable that the UCR for SigL exhibits large non-zero values in some intervals. An interesting observation from Panel D is that SL\_prod and SL\_min outperform ASigL in both the BA and WS networks.

In conclusion, we find that SigL underperforms compared to the other three methods in both sparse and non-sparse signal scenarios, exhibiting the lowest MCCa and non-zero UCR values, particularly at shorter lengths of $\Delta$ and in situations with high noise variance (as highlighted in Panel B). ASigL slightly outperforms SL\_prod, and SL\_min demonstrates robustness against noise variance in both sparse and non-sparse signal cases. For dense signals, the performance of SigL varies depending on the network type. While it performs best on BA networks, it is generally on par with SL\_prod and SL\_min, each excelling under specific conditions. Notably, SigL consistently shows larger non-zero UCR values in some intervals. On the other hand, the UCR values for ASigL, SL\_prod, and SL\_min are very close to zero across all cases, which is a desirable property. Between SL\_prod and SL\_min, we recommend using SL\_prod due to its greater robustness under specific conditions.

\subsection{\label{sec4} Kuramoto model in synchronization problem}

For problems introduced in Eq. (\ref{eq0}), we use the Kuramoto model  (\cite{r3}; \cite{r13}; Wu, {\it et al.}, 2012) to illustrate the reconstruction of the network in a complex system. This model has the following governing equation:
\begin{equation}\label{eq67}
\frac{d\theta_i}{dt}=\omega_i +c\sum_{j=1}^N a_{ij} sin(\theta_j -\theta_i),
\end{equation}
$i=1,\cdots, N$, where the system is composed of $N$ oscillators with phase $\theta_i$ and coupling strength $c$, each of the oscillators has its own intrinsic natural frequency $\omega_i$, $a_{ij}$ is the adjacency matrix of a give network and is need to be estimated in network reconstruction. Using the same framework of reference (\cite{r39}), the Euler method can be employed to generate time series with an equal time step $h$.  Let $Y_i=(y_{i1}, \cdots, y_{iL})'$, $y_{it}=[\theta_i(t+h)-\theta_i(t)]/h$, ${\bf X}_i=(\phi_{ij}(t))$ is a $L\times N$ matrix, with elements $\phi_{ij}(t)=c \times sin(\theta_j (t)-\theta_i(t))$
for $t=1,\cdots,L$ and $j=1,\cdots,N$, $\tilde \beta_i=(a_{i1}, \cdots, a_{iN})'$, then reconstruction model can be rewritten as
\begin{equation}\label{eq69}
Y_i=\omega_i {\bf 1}_L+{\bf X}_i \tilde \beta_i+\epsilon_i, 
\end{equation}
where ${\bf 1}_L$ denote a $L\times 1$ vector with all element 1.  

\begin{figure*}
\centering{\includegraphics[width=0.8\textwidth]{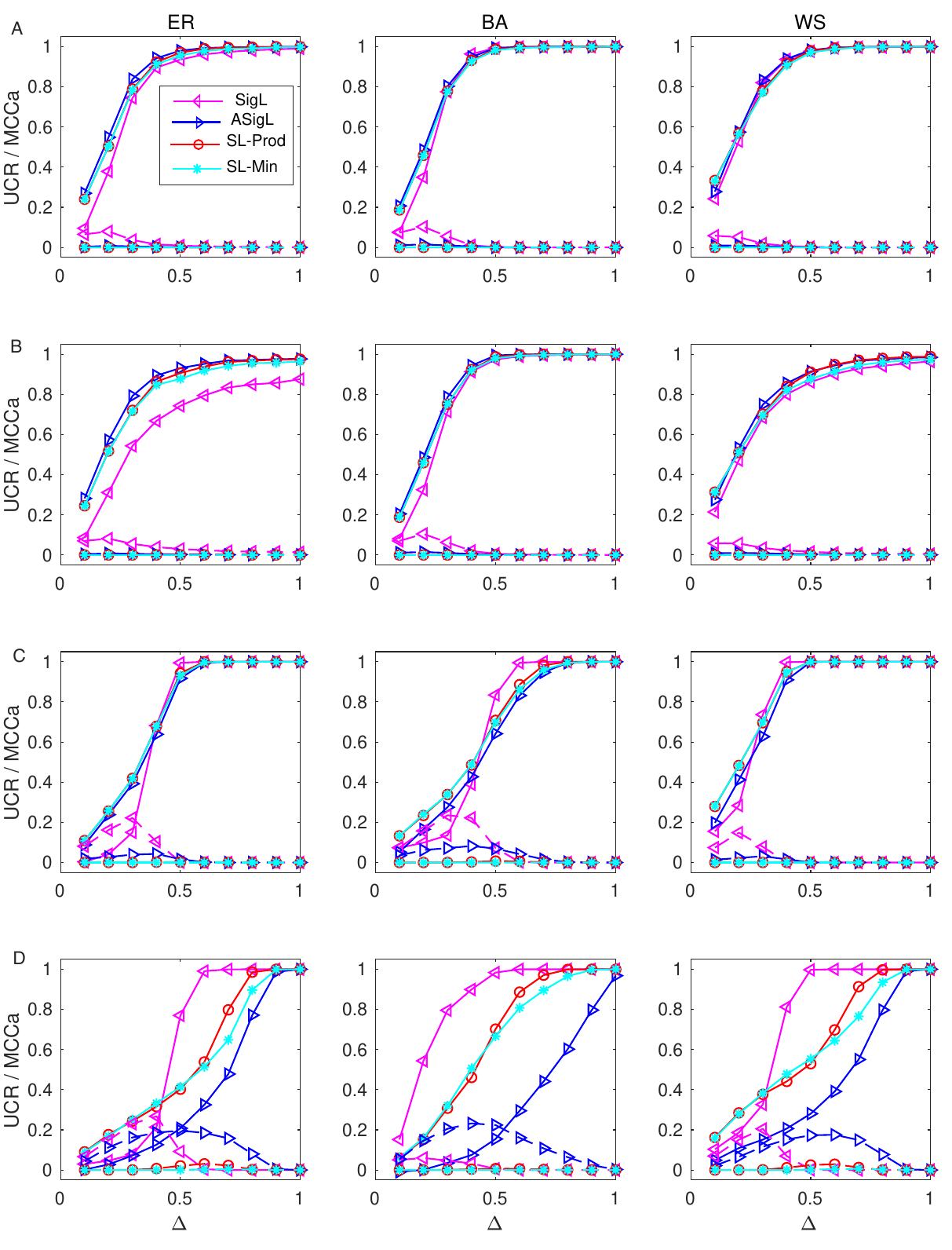}}
\caption{Accuracy measures MCCa and UCR in the reconstruction vs. $\Delta=L/N$, for Kuramoto model attained by  method of SigL, ASigL, SL\_prod and SL\_min in three kinds of network. The panel A contains the measures from the Kuramoto model in the network with average 6, while panel B contains the results for the case with noise of $\sigma^2 = 0.3$, and the network with average 6. The  panel C gives the results for non-sparse network with average degree 20 and without noise. The panel D refers to the results of dense signal for PDG game with average degree 70 and variance of noise $\sigma^2 = 0.1$. Three columns correspond to the results
based on Erd\"os-R\'enyi (ER) random networks,  Barab\'asi-Albert (BA) scale-free networks and  small world (WS) network, respectively. The network size $N = 100$ and coupling
c=10 for all cases. Each point is averaged over 20 simulations. }
\label{fig:6}    
\end{figure*}

Now we study the performance of the four methods in terms of Kuramoto model in ER, WS and BA networks with $N = 100$ and coupling strength $c=10$. The results are listed in Fig. \ref{fig:6}, where four panels A-D have the same network structures and noise perturbations as that in Section \ref{sec3}, but employing Kuramoto model as dynamics.  The overall results are similar to the evolutionary-game-based dynamical model.
In panel A and B, ASigL is slightly better than  SL\_prod and SL\_min. However the situation is reversed in panel C and D with methods of  SL\_prod and SL\_min being superior to ASigL especially in panel D.  Different with case of PD game in Section \ref{sec3}, the ASigL have non-zero UCR in panel C and D, while  UCR from SL\_prod and SL\_min remain the very small values (close to zeros). 
In Panel A, B and C, the performance of  SL\_prod and SL\_min almost coincide, but in Panel D, SL\_prod outperforms SL\_min for larger $\Delta$. Performance of four methods in panel D are very similar to Fig. \ref{fig:5}. The signal lasso have largest values of MCCa for BA network, while in ER and WS network, SigL have smaller MCCa when $\Delta <0.5$, but larger MCCa when $\Delta>0.5$. However we find that UCR for SigL have large non-zero values in some intervals.

Summing up the simulation results from linear regression, the evolutionary-game-based dynamical model and the Kuramoto dynamic model, the SL\_prod emerges as being on par with the ASigL for both sparse and non-sparse signals. The ASigL displays enhanced robustness against noise variance and can be used for noisy data. Yet, in scenarios involving dense signals, the SigL has some merit for larger $\Delta$, however its UCR will present some non-zeros values in some intervals. For dense signals, ASigL performs worse, while the performance of both SL\_prod and SL\_min remains stable and exhibits good overall performance. Moreover, the SL\_prod slightly outperforms its counterpart with the minimum penalty (SL\_min) when the design matrix reveals a zeros column. A significant advantage of both SL\_prod and SL\_min is the absence of a requirement to tune parameters, translating into substantial computational savings compared to the ASigL. We have calculate the time cost for four methods in different scenarios(see Table \ref{tb:A2} in Appendix), which show SL\_prod and SL\_min are faster and can be recommended for application in some complex situations.  

\section{Real Example: A human behavioral data}\label{sec5}

\begin{table}
    \caption{Payoff matrices used in the experiments. The payoff matrix of prisoner's dilemma game with punishment option (treatment I) and the standard prisoner's dilemma (treatment II and treatment III) are shown in left panel and right panel, respectively.}
    \label{tab:real2}
    \hrule
    \begin{minipage}{.5\linewidth}
      \centering
      \setlength{\tabcolsep}{6mm}{
        \begin{tabular}{cccc}
            & C & D & P \\
            \hline
            C & 2 & -2 & -5 \\
            D & 4 & 0 & -3  \\
            P & 2 & -2 & -5  \\
            \hline
            \end{tabular}}
    \end{minipage}%
    \begin{minipage}{.5\linewidth}
      \centering
       \setlength{\tabcolsep}{4mm}{
        \begin{tabular}{ccc} 
            ~ & C & D  \\
            \hline
            & & \\
            C & 4 & -2  \\
            D & 6 & 0   \\
            \hline
            \end{tabular}}
    \end{minipage} \\
    \hrule
\end{table}

In this section, we present the results of social network reconstruction using a real data from human behavioral experiment (\cite{r19}), where the purpose is to study the impact of the punishment on network reciprocity. A total of 135 participants from Yunnan University of Finance and Economics and Tianjin University of Finance Economics took part in the experiments and three trials are separately designed and carried out. In treatment I, 35 participants from Tianjin university of finance and economics played iterative prisoner's dilemma game with punishment on the static ring network with four neighbors, where there are 35 nodes and 140 links. For the other two treatments, 100 participants from Yunnan university of finance and economics were invited to participant in the iterated prisoner's dilemma experiment. The treatment II was implemented on the homogeneous random network with degree of 4, where there are 50 nodes and 200 links. While in treatment III, each player was placed on the heterogeneous random network, in which the degree of half nodes is 3 and the degree of other half nodes is 5, there are 50 nodes and 200 links. 
The network structure for three treatments are shown in \cite{r39} and the payoff matrices are listed in Table \ref{tab:real2}. In each round of the treatment, each player played with its direct neighbors to gain their payoff and updated its strategy to optimize its future payoff. The number of interactions in each session was set to 50, and the number of interactions was undisclosed until the sessions ended. In order to solve the network reconstruction problem in these examples, we thus recorded all the strategies and payoffs generated in the experiment and these information is the available data base for the reconstruction. This data set has been used by \cite{r39} and \cite{r40} to assess the network reconstruction of signal lasso and adaptive signal lasso. 

\begin{figure} 
\centering{\includegraphics[width=0.8\textwidth]{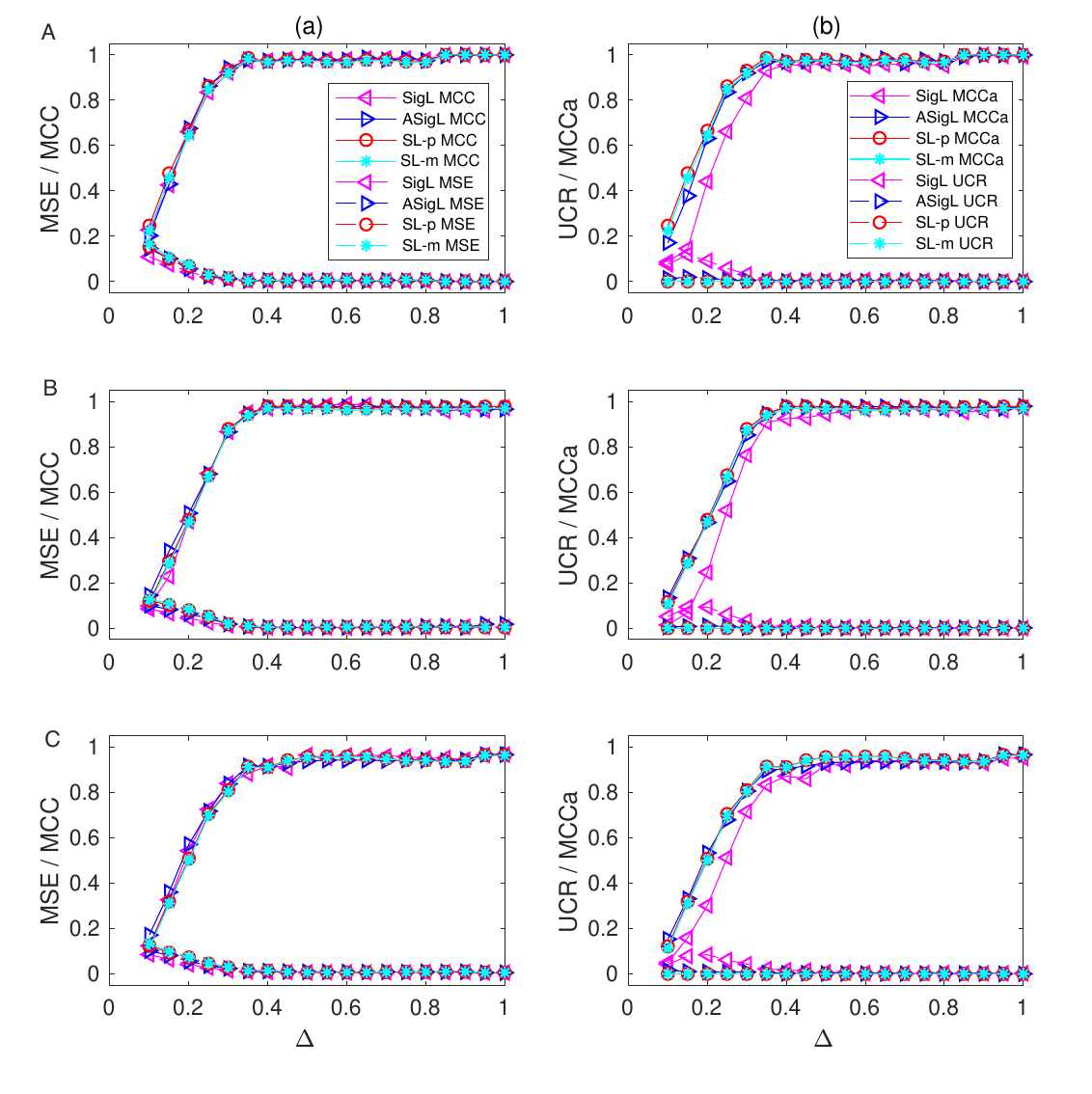}}
\caption{Accuracy in the reconstruction vs. $\Delta=L/N$, for three real trials using  method of SigL and ASigL.  Panel A refers to the results of the experimental ring network, where Panels A(a) is plot of MCC and MSE vs. $\Delta$, and A(b) refer to the MCCa and UCR criterion. There are 35 nodes and 140 links and the degree of each node is 4. Panel B refers to the results of the experimental homogeneous random network. There are 50 nodes and 200 links and the degree of each node is 4. Panel C refers to the results of the experimental heterogeneous random network. There are 50 nodes and 200 links and the average degree of each node is 4.}
\label{fig:7}    
\end{figure}

We again use this data set to compare the accuracy of reconstruction of network using different methods. The results are summarized in Fig.\ref{fig:7}, where upper panel A is about  results based on treatment I,  middle panel B is the  results based on treatment II, and  bottom panel C is for treatment III. 
Left column in the Fig.\ref{fig:7} give the results using measures of MCC and MSE.  We find that SL\_prod performs best in Fig.\ref{fig:7}A(a), while ASigL performs better slightly than other three methods  in Fig.\ref{fig:7}B(a) and Fig.\ref{fig:7}C(a) based on MCC. MSE of SigL have smaller values than other three methods for small $\Delta$ because other three method can shrink parameters to either 0 or 1 thus will present large deviation once it make a wrong selection.  

The right column  in Fig.\ref{fig:7} consider the influence of unclassified portion and thus use new measures of MCCa and UCR. The results in Fig.\ref{fig:7}A(b) show that SL\_prod performs best, and both SL\_prod and SL\_min have larger values of MCCa than that of ASigL and SigL. For treatment II and III as shown in Fig.\ref{fig:7}B(b) and Fig.\ref{fig:7}C(b), three methods, except for SigL, almost have the same values of MCCa. These three methods have approximately zeros values of UCR, where both SL\_prod and SL\_min almost exactly equal to zeros even for small $\Delta$. The signal lasso performs worse in these two treatments. 

\section{Conclusions}

In this paper we propose signal lasso estimation with two non-convex penalties, and make an extensive comparison with other shrinking methods. 
The results show that the performance of signal-lasso-type methods depends on the variance of noise and the density of network connectivity, but are robust against distribution of the data. 
We also discuss the performance of signal-lasso-type method in cases of non-sparse or dense signal, which usually appear in world trade web (WTW) (see \cite{r46}; \cite{r40}). 

In the simulations using linear regression models for both sparse and non-sparse signals, we observe that SL\_prod and SL\_min demonstrate superior or at least comparable performance to other methods in accurately recovering the network structure under conditions of low noise variance. However, with a large noise variance, the adaptive signal lasso performs best, especially when $n<p$, suggesting its robustness against noise perturbations in discerning complex network topologies.
In scenarios with dense signals where $n<p$, the signal lasso fares well for minor variances. Yet, as variance increase, SL\_prod and SL\_min outperform the rest.

Data sourced from evolutionary-game-based dynamical models and the Kuramoto model exhibit greater complexity compared to linear regression models. This data, often sparse within the design matrix, sometimes shows conditional illness or zero-columns. Our analysis under these conditions reveals that the performance of the four methods is critically dependent on network connectivity density. For sparse signals, the adaptive signal lasso (ASigL) is the most proficient, followed closely by SL\_prod and SL\_min. The signal lasso (SigL) lags, particularly under high variance conditions. In this scenario, ASigL slightly surpasses SL\_prod, while SL\_prod demonstrates superiority over SL\_min in the WS network. As network connectivity density increases, ASigL's advantages diminish, whereas SigL's effectiveness improves but is unstable. The results indicate that the performance of SigL and ASigL is unstable in dense networks, while the robust performance of both SL\_prod and SL\_min remains consistent. A notable advantage of both SL\_prod and SL\_min is their ability to bypass the need for selecting the tuning parameter 
$\lambda$, thereby speeding up computational times relative to ASigL (see Table \ref{tb:A2} in Appendix for detailed calculations). Furthermore, SL\_prod and SL\_min consistently shrink the parameter directly to 0 or 1, effectively eliminating any undetermined segments in classification challenges. 

In conclusion, our results show that signal lasso with non-convex penalties is effective and fast  in estimating signal parameters in linear regression model. 
We recommend using SL\_prod or ASigL in sparse networks, with ASigL being the preferred choice when the data exhibit high variance or significant noise. For non-sparse or dense networks, SL\_prod and SL\_min are the preferred methods. Between SL\_prod and SL\_min, we suggest using SL\_prod, as it is more robust against certain special conditions that may arise in the design matrix.

\section*{Acknowledgement}
This work was supported by National Natural Science Foundation of China (NSFC, Grant No. 12271471, 11931015),   National Social Science Foundation of China (Grants Nos. 22\&ZD158 and 22VRCO49) and Project of "XingDianYingCai" Plan of Yunnan Province to L. S. and JSPS Postdoctoral Fellowship Program for Foreign Researchers (No. P21374) and an accompanying Grant-in-Aid for Scientific Research to C.\,S.

\section*{Appendix}
\renewcommand{\theequation}{A.\arabic{equation}}

\subsection{The results for zero column in design matrix of linear regression model}
In Table \ref{tb:A1}, we list the simulation results in linear regression model with design matrix $X$ have 2 zero columns for nine methods: Lasso, Adaptive Lasso, SCAD, MCP, Elastic Net, SigL, ASigL, SL\_prod and SL\-min in linear regression models. 
All of the results are averaged over 500 independent realizations, where n is the sample size, $p$ is the number of explanatory variables, $p_1$ is the number of signals (number of $\beta=1$). 
In this special case, it is clear that  SL\_prod has better performance than that of  SL\_min. 

\begin{table}[h]\footnotesize
	\begin{center}
		\caption{Simulation results in linear regression model with design matrix $X$ have 2 zero columns for nine methods}
		\label{tb:A1} 
\setlength{\tabcolsep}{3mm}{
	\begin{tabular}{ccccccccc}	
			\hline
Method	&  \multicolumn{3}{c}{MSE/UCR/MCC/MCCa }  \\
\hline
$(n, p, p_1, \sigma)$ 		&$(50, 150, 6, 0.4)$ & $(50, 150, 6, 1)$  & $(50, 150, 6, 2)$\\	
\hline
Lasso &0.0020/0.049/0.960/0.435  &  0.011/0.173/0.716/0.036 & 0.0423/0.265/0.208/-0.090  \\
A-lasso &0.0024/0.071/0.960/0.346 & 0.0146/0.193/0.664/-0.001 &0.0464/0.269/0.291/-0.081 \\
SCAD & 0.0002/0.004/1.000/0.934  &  0.0059/0.036/0.799/0.394 &0.0312/0.036/0.067/-0.003  \\
MCP & 0.0008/0.005/0.980/0.911  &  0.0065/0.030/0.811/0.426 & 0.0332/0.025/0.110/0.031\\
ElasticNet&0.0019/0.034/0.780/0.365  & 0.0113/0.062/0.316/0.061 &0.0316/0.038/0.048/-0.002 \\ 
SigL &0.0009/0.026/1.000/0.640 & 0.0069/0.098/0.950/0.188& 0.0281/0.128/0.475/0.021\\
ASigL &0.0103/0.005/0.862/0.805  &  0.0156/0.004/0.779/0.736 &0.0314/0.006/0.525/0.479\\
SL\-prod &   0.0000/0.000/0.999/1.000  & 0.0058/0.000/0.936/0.936 & 0.0891/0.000/0.386/0.386\\
SL\-min &   0.0133/0.000/0.859/0.859  & 0.0191/0.000/0.806/0.806& 0.1024/0.000/0.355/0.355\\    

\hline
		\end{tabular}
		}
	\end{center}
\end{table}

\subsection{Computational complexity measures (CPU)}
In Table \ref{tb:A2} , we calculated the CPU time (seconds) for computing estimations of parameters based on  SigL, ASigL, SL\_prod and SL\_min methods, where the initial value for computation is lasso estimators for all methods. It is clear that SL\_prod and SL\_min methods
have obvious advantage in computing speed.

\begin{table}[h]\footnotesize
\centering
\caption{Comparison of computational complexity measures (CPU) for four methods}
\label{tb:A2} 
\begin{tabular}{ccccccc}
\toprule
\multicolumn{7}{c}{Linear Dynamic ($\sigma=2$, Correlated-Gaussian)} \\
\hline
$(n, p, p_1)$ & (150,30,6) & (150,150,6) & (150,150,60) & (200,30,6) & (200,150,6) & (200,150,60) \\
\hline
SigL & 0.389935 & 16.330999 & 7.049003 & 0.466000 & 11.325395 & 6.643295 \\
ASigL & 0.130004 & 0.457998 & 0.503997 & 0.157996 & 0.556997 & 0.644602 \\
SL\_prod & 0.000000 & 0.000000 & 0.000999 & 0.000000 & 0.001004 & 0.000995 \\
SL\_min & 0.001000 & 0.001003 & 0.001001 & 0.000000 & 0.001000 & 0.001004 \\
\hline
\multicolumn{7}{c}{PDG Dynamic with $N=100, \Delta=0.5, b=1.02$} \\
\hline
 & \multicolumn{3}{c}{k=6} & \multicolumn{3}{c}{k=60} \\
\hline
 & BA & ER & SW & BA & ER & SW \\
\hline
SigL & 1.138426 & 0.848676 & 0.569877 & 2.514897 & 3.014704 & 2.714757 \\
ASigL & 0.044129 & 0.039120 & 0.039770 & 0.041260 & 0.039300 & 0.039610 \\
SL\_prod & 0.000441 & 0.000502 & 0.000363 & 0.000374 & 0.000452 & 0.000372 \\
SL\_min & 0.000438 & 0.000503 & 0.000361 & 0.000444 & 0.000371 & 0.000453 \\
\hline
\multicolumn{7}{c}{Kuramoto Dynamic with $N=100,\Delta=0.5, \lambda=15$} \\
\hline
 & \multicolumn{3}{c}{k=6} & \multicolumn{3}{c}{k=60} \\
\hline
 & BA & ER & SW & BA & ER & SW \\
\hline
SigL & 2.433161 & 4.215170 & 4.553449 & 3.071393 & 4.965631 & 3.801474 \\
ASigL & 0.118870 & 0.260713 &	0.306684 &	0.202545 &	0.304234 &	0.191870 \\
SL\_prod& 0.000580 & 0.001184 &	0.001225 &	0.001444 &	0.003402 &	0.001280 \\
SL\_min & 0.000423 & 0.000365 &	0.000367 &	0.000375 &	0.000380 &	0.000378 \\
\hline
\end{tabular}
\end{table}

\subsection{Program of algorithm and real dataset}
The computation of methods such as lasso,
adaptive lasso, SCAD, MCP, and ElasticNet can be found in
R-software and has been widely used in scientific research. The R-code using
the coordinate descent method for the SigL, ASigL, SL\_prod and SL\_min, as well as the real dataset are given on GitHub https://github.com/shilei65/adaptive-signal-lasso-code.git.

\end{document}